\documentclass[12pt]{amsart}
\usepackage[final]{graphicx}
\usepackage{cite}
\usepackage{booktabs}
\usepackage{enumerate}
\usepackage[tight]{subfigure}
\usepackage{url}
\usepackage[bookmarks=true, bookmarksopen=true,%
    bookmarksdepth=3,bookmarksopenlevel=2,%
    colorlinks=true,%
    linkcolor=blue,%
    citecolor=blue,%
    filecolor=blue,%
    menucolor=blue,%
    urlcolor=blue]{hyperref}
\usepackage[margin=1in]{geometry}
\usepackage[textsize=small]{todonotes}

\newcommand{\RR}{\mathbb R}
\newcommand{\CC}{\mathbb C}

\theoremstyle{plain}
\newtheorem{proposition}{Proposition}
\newtheorem{theorem}[proposition]{Theorem}
\newtheorem{lemma}[proposition]{Lemma}
\newtheorem{corollary}[proposition]{Corollary}

\newtheorem{conjecture}[proposition]{Conjecture}

\theoremstyle{definition}
\newtheorem{definition}[proposition]{Definition}

\newtheorem{question}[proposition]{Question}

\theoremstyle{remark}

\newtheorem*{remark}{Remark}
\def\mathcenter#1{%
  \vcenter{\hbox{#1}}%
}

\def\mfig#1{
        \mathcenter{\includegraphics{#1}}
}

\def\mfigb#1{
        \mathcenter{\includegraphics[trim=-1 -1 -1 -1]{#1}}
}

\graphicspath{{draws/}{figs/}}

\begin{document}
\title{From dominoes to hexagons}
\author[Thurston]{Dylan~P.~Thurston}
\address{Department of Mathematics\\
         Indiana University, Bloomington\\
         831 E. Third St.,
         Bloomington, Indiana 47405\\
         USA}
\email{dpthurst@indiana.edu}
\date{September 12, 2016}

\begin{abstract}
  There is a natural generalization of domino tilings to tilings of a
  polygon by hexagons, or, dually, configurations of oriented curves
  that meet in triples.  We show exactly when two such tilings can be
  connected by a series of moves analogous to the domino flip move.
  The triple diagrams that result have connections to Legendrian
  knots, cluster algebras, and planar algebras.
\end{abstract}

\maketitle

\section{Introduction}
\label{sec:intro}
\begin{figure}[htbp]
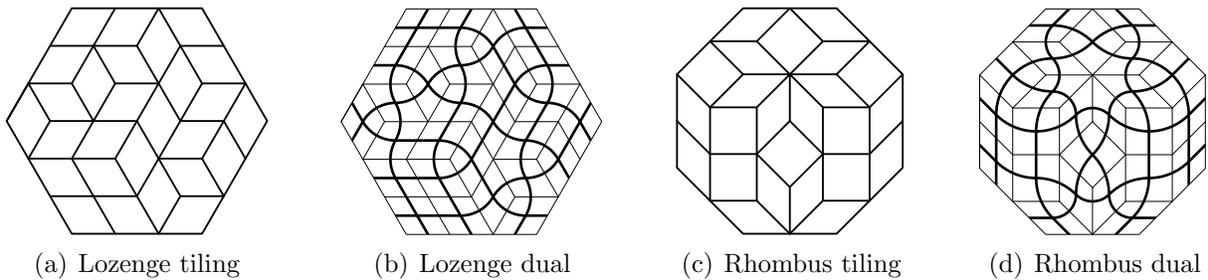

  \centerline{
  \subfigure[Lozenge tiling]
    {\label{fig:lozenge}\includegraphics{lozenge-0}}
  \qquad
  \subfigure[Lozenge dual]
    {\label{fig:lozenge-dual}\includegraphics{lozenge-1}}
  \qquad
  \subfigure[Rhombus tiling]
    {\label{fig:rhombus}\includegraphics{lozenge-10}}
  \qquad
  \subfigure[Rhombus dual]
    {\label{fig:rhombus-dual}\includegraphics{lozenge-11}}
  }
  \caption{Tilings by lozenges or more general rhombuses and their duals.}
  \label{fig:classical}
\end{figure}

The study of tilings of a planar region by lozenges, as in
Figure~\ref{fig:lozenge}, or more generally by rhombuses, as in
Figure~\ref{fig:rhombus}, has a long history in combinatorics.
Interesting questions include deciding when a region can be tiled,
connectivity of the space of tilings under the basic move
$\mfigb{lozenge-20}\leftrightarrow\mfigb{lozenge-21}$, counting the
number of tilings, and behavior of a random tiling.

One useful tool for studying these tilings is the dual picture, as in
Figures~\ref{fig:lozenge-dual}. It consists
of replacing each rhombus in the tiling by a cross of two strands
connecting opposite, parallel sides.  If we trace a strand through the
tiling it comes out on a parallel face on the opposite side of the
region, independently of the tiling of the region. Lozenge tilings
give patterns of strands with three families of parallel
(non-intersecting strands).

If we drop the restriction to three families of parallel strands, as
in Figure~\ref{fig:rhombus-dual}, and take any set of strands
connecting boundary points
which are both \emph{generic} (only two
strands intersect at a time) and \emph{minimal} (no strand intersects
itself and no two strands intersect more than once), then there is a
corresponding dual tiling by rhombuses, as in
Figure~\ref{fig:rhombus}. (This was proved in great generality by
Kenyon and Schlenker \cite{KS05:RhombicEmb}.)

These rhombus tilings correspond to words in the symmetric group, and
any two rhombus tilings of a region can be related by Reidemeister III
moves on the strands diagram, which in terms of rhombus tilings
replaces a hexagon tiled with 3 rhombus tiles with another tiling of
the same hexagon by the same 3 rhombi. This has been proved by several
people, including by Curtis, Ingerman, and Morrow in their work on
resistor networks \cite{CIM98:CircNetworks}.

\begin{figure}[htbp]
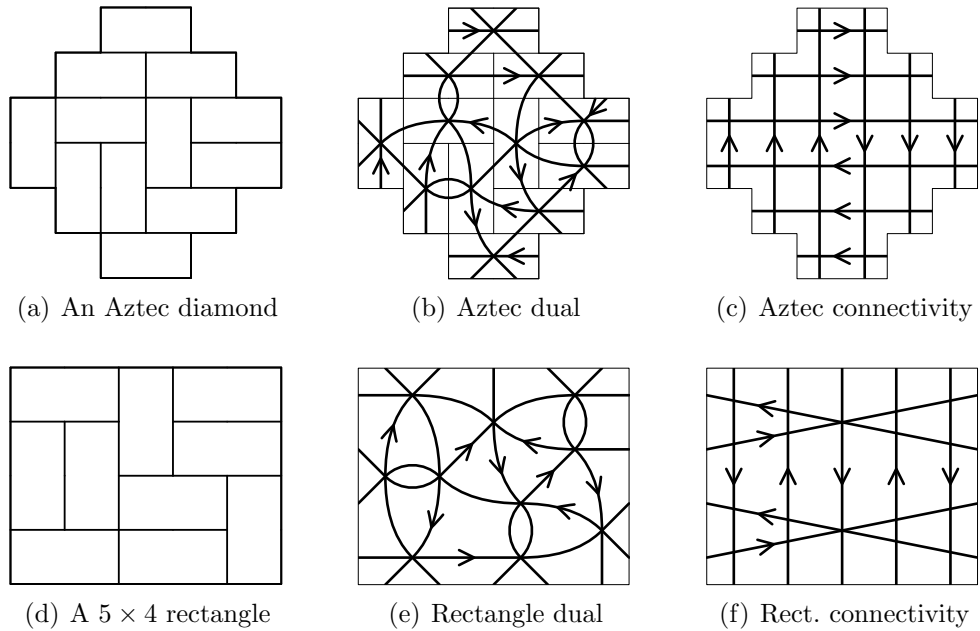

  \centerline{
    \subfigure[An Aztec diamond]
      {\label{fig:aztec}\includegraphics{domino-0}}
    \qquad
    \subfigure[Aztec dual]
      {\label{fig:aztec-dual}\includegraphics{domino-101}}
    \qquad
    \subfigure[Aztec connectivity]
      {\label{fig:aztec-connect}\includegraphics{domino-103}}
  }
  \vskip 5pt
  \centerline{
    \subfigure[A $5\times4$ rectangle]
      {\label{fig:rectangle}\includegraphics{domino-10}}
    \qquad
    \subfigure[Rectangle dual]
      {\label{fig:rectangle-dual}\includegraphics{domino-12}}
    \qquad
    \subfigure[Rect.~connectivity]
      {\label{fig:rectangle-connect}\includegraphics{domino-14}}
  }
  \caption{Domino tilings and their duals, along with the connectivity
    of the strands}
  \label{fig:domino}
\end{figure}

In this paper, we will look at a related situation, in which tilings
by dominoes
(Figure~\ref{fig:aztec}) are given a dual picture
(Figure~\ref{fig:aztec-dual}), in which we replace each domino by an
``asterisk'' of three strands which are not generic, but rather meet
in a triple point.  When we perform the basic \emph{domino flip}
\[
\mfigb{domino-20}\leftrightarrow\mfigb{domino-21}
\]
on domino tilings, the dual strands perform what we will call a
\emph{$2\leftrightarrow2$ move}
\[
\mfigb{domino-25}\leftrightarrow\mfigb{domino-26}.
\]
In particular, the connectivity of the strands as we trace them across
the diagram is unchanged.  Since the space of domino tilings is
connected under the domino flip~\cite{Thurston90:ConwaysTiling}, the
connectivity of the strands depends only on the shape of the region
and not on the particular tiling (Figure~\ref{fig:aztec-connect}).

Conversely, it follows from results in this paper that this is a
faithful representation: every set of immersed arcs with the same
connectivity as a domino tiling, meeting only in triple points, and
with a minimal number of triple points corresponds uniquely to a domino
tiling.

As in the case of rhombus tilings, strands always connect parallel
sides.  Unlike in the rhombus case, strands may intersect each other
arbitrarily many times.  For the special case of the tiling of the
Aztec diamond (Figures~\ref{fig:aztec}--\ref{fig:aztec-connect}),
strands connecting parallel sides do not cross in the net connectivity
across the diagram.  But this is not true in general, as you can see
in Figures~\ref{fig:rectangle}--\ref{fig:rectangle-connect}.

But there is something we can say about the crossings of strands.
Consider the complementary regions to the strands.  Because they meet
with six around a vertex, they may be checkerboard colored with two
colors, black and white.  Orient the strand segments clockwise around
the black regions and counterclockwise around the white regions.
Because 6 is congruent to 2 modulo 4, the orientations are consistent
when we follow a strand through a triple crossing:
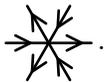
\[
\mfigb{arcs-0}.
\]
See Figure~\ref{fig:rectangle-dual} for the orientations of the strands
in the dual to the domino tiling in Figure~\ref{fig:rectangle}.
Another way to characterize these orientations is to consider a
generic curve inside the diagram; the strands that it crosses will
alternate crossing left-to-right and right-to-left.

The appropriate analogue of the minimal condition for the diagrams
dual to rhombus tilings is that strands connecting parallel sides
oriented in the same direction never cross each other.

These pictures suggest a natural generalization.

\begin{definition}
  A \emph{triple diagram} is a collection of oriented one-manifolds,
  possibly with boundary, mapped smoothly into the disk.  The image of
  a connected component is a \emph{strand}; it is either an \emph{arc}
  (the image of an interval) or a \emph{loop} (the image of a circle).
  The maps are required to satisfy:
  \begin{itemize}
  \item Three strands cross at each point of intersection;
  \item the endpoints of arcs are distinct points on the boundary of
    the disk, and no other points map to the boundary; and
  \item the orientations on the strands induce consistent orientations
    on the complementary regions.
  \end{itemize}
  Triple diagrams are considered up to homotopy among such diagrams.
  This makes them essentially combinatorial objects, 6-valent graphs
  with some extra structure.
\end{definition}

%\begin{definition}
%  A \emph{triple diagram} is an oriented 1-manifold properly
%  immersed in the disk so that three strands cross at every
%  intersection point, the orientations on the strands induce a
%  consistent orientation on each complementary region, and only one
%  strand meets the boundary at any point.  The connected components of
%  the 1-manifold are called the \emph{strands}.  They are of two
%  types: the \emph{arcs} which meet the boundary (topologically
%  intervals), and the \emph{loops} which do not (topologically
%  circles).
%\end{definition}

\begin{definition}\label{def:connected}
  A \emph{connected} triple diagram~$D$ is a diagram in which the
  image of the immersed curves together with the boundary of the disk
  is connected.  Equivalently, it is a diagram in which each
  complementary region to the image is a disk.  A \emph{disjoint
    component} of a triple diagram is a connected component of the
  image of the immersed curves which does not meet the boundary of the
  disk.  A \emph{simple loop} is a loop which goes through no simple
  crossings and bounds an empty disk.
\end{definition}

A diagram in which all strands are arcs is automatically connected.

If the diagram is connected, the condition on orientations amounts to
requiring that the orientations alternate in and out around each
triple point, and is automatically satisfied if we start anywhere and
start assigning alternating orientations around vertices.  Since all
the complementary regions of a connected diagram are disks, we can
construct the dual, which is a topological tiling by hexagons.  Note
that a domino tiling can be turned into a tiling by hexagons by adding
a vertex in the middle of the long edges of each domino.  The triple
crossing diagram described above is dual to this tiling by
hexagons.

As before, the connectivity of the arcs gives a matching of the
vertices on the boundary of the disk.
A natural question to ask is which matchings of the boundary vertices
are achievable by a triple diagram.  There is one immediate
restriction. The orientations on the strands alternate in and out as
we go around the boundary of the disk, and every inward-pointing end
gets matched with an outward-pointing end.  The matching is a bijection
between the in- and out-endpoints.

\begin{theorem}\label{thm:generate}
  In a disk with $2n$ endpoints on the boundary, all $n!$ pairings of
  in-endpoints with out-endpoints are achievable by some triple point
  diagram without closed strands.
\end{theorem}

The proof is in Section~\ref{sec:generating}.

If we think about these diagrams by analogy with group theory, this
theorem assures us that triple crossings \emph{generate} admissible
pairings.  The next question is \emph{relations}: What moves on
diagrams of triple points allow us to pass between two diagrams which
induce the same matching?  One natural move is the $2\leftrightarrow2$ move above.
This move suffices if the total number of triple points is as small as
possible.

\begin{definition}
  A \emph{minimal} triple diagram is a connected diagram with no more
  triple points than any other triple diagram inducing the same
  pairing on the boundary.
\end{definition}

(The condition of connectivity merely rules out small loops with no crossings
disconnected from the rest of the diagram.)

\begin{theorem}\label{thm:domino-flip}
  Any two minimal triple point diagrams with the same matching on the
  endpoints can be related by a sequence of $2\leftrightarrow2$ moves.
\end{theorem}

The restriction to minimal diagrams is necessary, since the $2\leftrightarrow2$
move, as per its name, keeps the total number of triple points
unchanged.  To drop this restriction, we need to introduce some moves
that change the number of triple points while preserving the matching.
These moves have no analogues in domino tilings, since the number of
triple crossings in a domino tiling is the number of dominoes, i.e.,
half the area of the region, which is invariant.

\begin{theorem}\label{thm:reduction}
  Any triple point diagram can be reduced to any minimal triple point
  diagram (with the same matching) by a sequence of the following
  moves:
  \begin{itemize}
  \item $2\leftrightarrow2$ moves $\mfigb{moves-0} \leftrightarrow
    \mfigb{moves-1}$;
  \item $1\to0$ moves $\mfigb{moves-11} \to
    \mfigb{moves-12}$; and
  \item Dropping a simple loop~$\mfigb{moves-70}$ with no
    crossings and an empty interior.
  \end{itemize}
\end{theorem}

Note that this version of the theorem is a little stronger than you
might expect: we never increase the number of triple crossings.

Theorems~\ref{thm:domino-flip} and~\ref{thm:reduction} are proved in
section~\ref{sec:reducing}. The result for domino tilings was
previously known \cite{Thurston90:ConwaysTiling,STCR95:SpacesDomino}.

\begin{remark}
  There is an alternate version of Theorem~\ref{thm:reduction} which
  stays within the space of connected diagrams: any connected diagram
  can be reduced to a minimal diagram by a sequence of connected
  diagrams related by
  \begin{itemize}
  \item $2\leftrightarrow2$ moves;
  \item $1\to0$ moves;
  \item Dropping move 1 $\mfigb{moves-20} \to
    \mfigb{moves-21}$; and
  \item Dropping move 2 $\mfigb{moves-30} \to
    \mfigb{moves-31}$.
  \end{itemize}
  Since this version stays entirely within connected diagrams, it can
  also be stated in the dual language of hexagons.  We will not prove
  this version here.  It can be deduced as a corollary of
  Theorem~\ref{thm:reduction} by analyzing the diagram just before
  applying a $1\to0$ move which disconnects the diagram.
\end{remark}

There is also an internal characterization of which triple point
diagrams are minimal.

\begin{theorem}\label{thm:minimal}
  A connected triple point diagram is minimal if and only if it has no
  strands which intersect themselves (\emph{monogons}) or pairs of
  strands which intersect at two points, $x$ and $y$, with both strands
  oriented from $x$ to $y$ (\emph{parallel bigons}).
\end{theorem}

The proof is in section~\ref{sec:minimal}.  Note that a disconnected
component is either a simple loop or has a monogon so the restriction
to connected diagrams is not a large restriction.

This is an analogue of the theorem that a generic immersion of arcs in
the plane has a minimal number of crossing points for the given
pairing of the boundary (and corresponds to a tiling by rhombuses) iff
there are no monogons or bigons.  Note that in the triple point
diagrams in Figure~\ref{fig:domino} you can see many bigons, but all
of these bigons are anti-parallel: if the intersections are at A and
B, one strand runs from A to B and the other strand runs from B to A.

\subsection*{History}
The first version of this paper was completed and posted on the
mathematics arXiv in 2004. The present version has a few minor
updates,
including an improved proof of
Theorem~\ref{thm:crossing-count} and motivation for
Conjecture~\ref{conj:rels-rels}.
After the first version of this paper was completed, I learned that
Alexander Postnikov independently arrived at many of the same results
in a somewhat different language. In particular, \cite[Theorem
12.7]{Postnikov06:TotalPositivity} is closely related to
Theorem~\ref{thm:domino-flip}, with his move~(M1) corresponding to our
$2 \leftrightarrow 2$ move.

\subsection*{Acknowledgements}

I would like to thank my advisor, Vaughan Jones, for suggesting this
problem, and for his patience in waiting for me to write it
up and then publish it. I am very pleased to be able to publish this
paper in a volume celebrating his birthday.

I would also like to thank Vladimir Arnol'd, Sergei Fomin, Andre
Henriques, Michael Polyak, Alexander Postnikov, Jim Propp, Chung-chieh
Shan, and Vladimir Turaev for many helpful discussions and comments.

This work was completed while the author was at UC Berkeley and at
Harvard University, and was partially supported by the National Science
Foundation under a Graduate Research Fellowship and Grant Numbers
DMS-0071550 and DMS-1507244.

\section{Generating matchings}
\label{sec:generating}

In order to prove Theorem~\ref{thm:generate}, we will explicitly
construct a triple diagram that achieves a given pairing of
in- and out-endpoints.  The triple diagrams we construct
turn out to be minimal diagrams.  We will use this construction in the
sections that follow: we will show that any two diagrams are related
by reducing both of them to a diagram of the form we construct here.

We are given a circle with $2n$ endpoints around the boundary,
alternately marked ``in'' and ``out'', and a pairing of the
in-endpoints with the out-endpoints.  Each pair of an in-endpoint and
an
out-endpoint divides the circle into two \emph{intervals}.  These
intervals are ordered by inclusion.  Pick a minimal
interval~$I$ with respect to this order.  Now start to construct the
triple diagram by running a \emph{boundary-parallel} strand~$S$
along~$I$, introducing a triple crossing for each pair of strands that
we cross over:
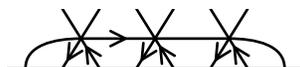
\[
\mfigb{arcs-10}
\]
Note that there will always be an even number (possibly zero) of
strands to cross
over because of the alternation of in- and out-endpoints.

Remove~$S$ from our original pairing and swap the pairs we crossed
over.  The in- and out-endpoints in the new pairing still alternate and
there are fewer strands. Continue by induction until all the endpoints
are paired up.

A diagram constructed in this way, for some sequence of choices of
minimal intervals, is called \emph{standard}. We will see later
(Corollary~\ref{cor:std-minimal}) that
standard diagrams are minimal.

\begin{remark}
  We could be a little more liberal in which intervals we allow, while
  still preserving minimality of the resulting triple crossing
  diagram.  We can distinguish between the two intervals arising from
  a given matching according to whether the interval runs \emph{clockwise} or
  \emph{counterclockwise} from the in-endpoint to the out-endpoint.  If, in
  the construction above, we
  always pick a clockwise interval that is minimal among all clockwise
  intervals, or a counterclockwise interval minimal among all
  counterclockwise intervals, the
  resulting triple diagram is still minimal.
\end{remark}

We can also count the number of crossings in a standard diagram using
only the connectivity information of the strands. For this purpose,
pick any linear ordering on the strand endpoints that is compatible
with the cyclic ordering. (That is, break the circle open at a
basepoint, and lay it out on a line.)

\begin{definition}
  A strand in a triple diagram from~$a$ to~$b$ (with a choice of
  basepoint) is \emph{right-moving}
  if $a < b$, and is \emph{left-moving} if $b < a$.
  Two strands $S$ in a triple diagram are \emph{linked} if their
  endpoints interleave around the circle. If the strands are linked,
  they are \emph{parallel} if they are both left-moving or both
  right-moving, and are \emph{anti-parallel} if one is left-moving and
  one is right-moving. Explicitly, two strands are parallel linked if
  their endpoints are in one of these two patterns:
  \[
\mfigb{arcs-20}\qquad\mfigb{arcs-21}
\]
  (Note that these definitions refer to the strands, but they
  depend only on the connectivity.)
\end{definition}

\begin{theorem}\label{thm:crossing-count}
  The number of triple crossings in a standard triple diagram~$D$ is
  equal to the number parallel linked of pairs of strands.
\end{theorem}

\begin{lemma}\label{lem:count-order}
  For a given strand connectivity, the number of parallel linked pairs
  of strands doesn't depend on where you cut the circle open (i.e.,
  which linear ordering compatible with the cyclic ordering you pick).
\end{lemma}

\begin{proof}
  Consider a strand~$S$ from $a$ to~$b$, with $b$ the largest value in
  the linear ordering. $S$ is necessarily right-moving. A strand that
  links~$S$ is right-moving if it has its out-endpoint between $a$
  and~$b$, and is left-moving if it has its in-endpoint between $a$
  and $b$. Thus $S$ is linked with an equal number of left-moving and
  right-moving strands, and we get the same number of parallel linked
  pairs of strands if we change the linear ordering so that $b$ is the
  smallest value. Other cases are symmetric.
\end{proof}

\begin{proof}
  We count the number of crossings in a standard diagram by induction
  on the number of strands.  Pick a minimal interval, say from~$a$
  to~$b$. By Lemma~\ref{lem:count-order}, we may suppose the interval
  does not contain the basepoint.  By symmetry, we may suppose
  $a < b$.  Now consider what happens when we run a boundary-parallel
  strand~$S$ from~$a$ to~$b$ crossing a pair of adjacent
  out-endpoint~$c$ and in-endpoint~$d$.  Let the two other strands
  be~$T$ and~$U$, so that $S$ is linked with a strand $T$ from another
  point~$e$ to~$c$ and a strand $U$ from~$d$ to another point~$f$.
  Let the diagram (with one fewer strand) above~$S$ be $D'$.  Consider
  the location of~$e$ and~$f$ relative to~$a$ and~$b$. (Note neither $e$
  nor~$f$ can be between $a$ and~$b$, since $T$ and~$U$ are linked
  with~$S$ by the assumption that $[a,b]$ was a minimal interval.)
  \begin{description}
  \item[Case 1. $e < a$, $f < a$] The linking between $S$ and~$T$ is
    parallel, the linking between $S$ and~$U$ is anti-parallel, and
    any linking created or destroyed between~$T$ and~$U$ is anti-parallel.
  \item[Case 2. $b < e$, $b < f$] The linking between~$S$ and~$T$ is
    anti-parallel, the linking between~$S$ and~$U$ is parallel, and
    any linking created or destroyed between~$T$ and~$U$ is anti-parallel.
  \item[Case 3. $f < a < b < e$] The linkings between~$S$ and~$T$
    and~$S$ and~$U$ are both anti-parallel, while $T$ and~$U$ are
    parallel linked in~$D$ but not in~$D'$.
  \item[Case 4. $e < a < b < f$] The linkings between~$S$ and~$T$
    and~$S$ and~$U$ are both parallel, while $T$ and~$U$ are parallel
    linked in~$D'$ but not in~$D$.
  \end{description}
  In each case, the number of parallel linked strands in~$D'$
  decreased by~$1$ relative to~$D$ as a result of this crossing.
\end{proof}

\section{Reducing diagrams}
\label{sec:reducing}

In order to prove Theorems~\ref{thm:domino-flip}
and~\ref{thm:reduction}, we will start with an arbitrary diagram and
reduce it to one of the form described in Section~\ref{sec:generating}
by straightening the strands one-by-one along the boundary.

The basic blocks in the proof are some relations which let
us slide one strand past another.

\begin{lemma}\label{lem:useful-moves}
  Each of the following triple-point diagrams can be related by
  sequences of $2\leftrightarrow2$ moves:
  \renewcommand\theenumi{\alph{enumi}}
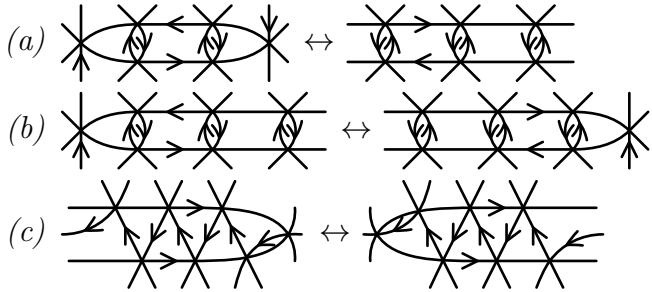
  \begin{enumerate}
  \item $\mfigb{domino-30} \leftrightarrow \mfigb{domino-35}$\label{it:move-1}
  \item $\mfigb{domino-40} \leftrightarrow \mfigb{domino-45}$\label{it:move-2}
  \item $\mfigb{domino-50} \leftrightarrow \mfigb{domino-55}$\label{it:move-3}
  \end{enumerate}
  In each diagram the central section can be repeated an arbitrary
  number of times.
\end{lemma}

\begin{proof}
  In each case, the two triple-point diagrams come from two different domino
  tilings of the same region:
  \renewcommand\theenumi{\alph{enumi}}
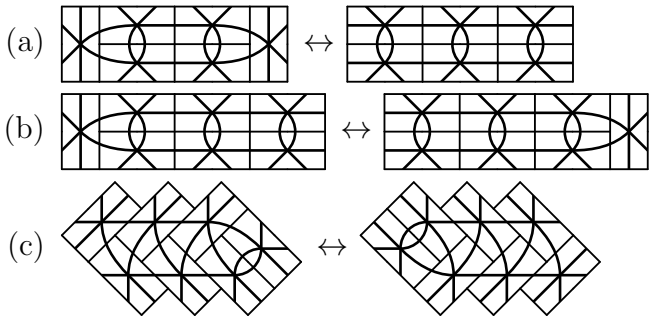
  \begin{enumerate}
  \item $\mfigb{domino-31} \leftrightarrow \mfigb{domino-36}$
  \item $\mfigb{domino-41} \leftrightarrow \mfigb{domino-46}$
  \item $\mfigb{domino-51} \leftrightarrow \mfigb{domino-56}$
  \end{enumerate}
  We use the fact that any two domino tilings of the same region can
  be related by a sequence of domino flips. (It is easy to find an
  explicit sequence of flips in these cases.)
\end{proof}

\begin{figure}[htbp]
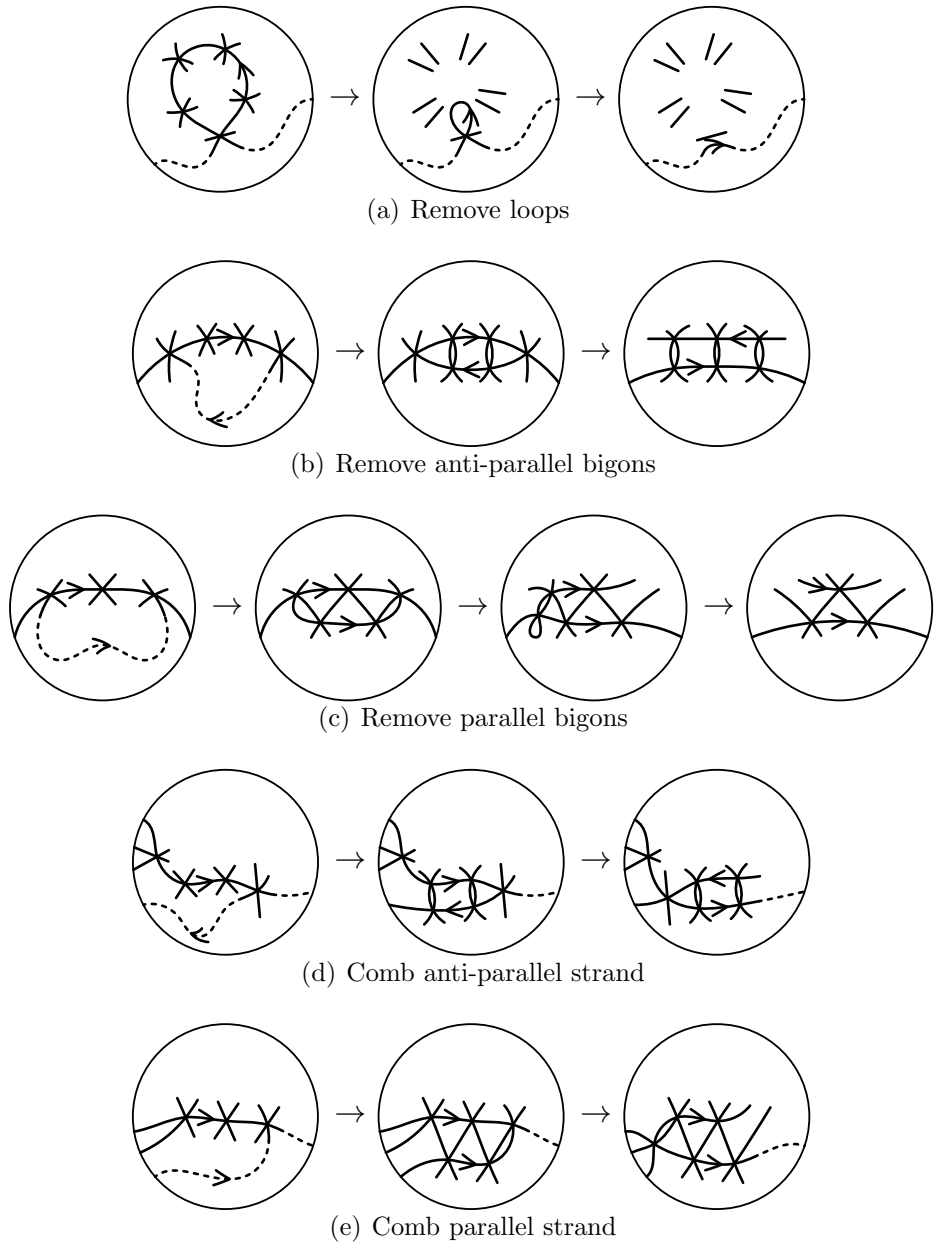

\centerline{
  \subfigure[Remove loops]
  {\label{fig:remove-loop}$\mfigb{simplify-10}\to\mfigb{simplify-11}\to\mfigb{simplify-12}$}
  }
\centerline{
  \subfigure[Remove anti-parallel bigons]
  {\label{fig:remove-anti-parallel-bigon}$\mfigb{simplify-0}\to\mfigb{simplify-1}\to\mfigb{simplify-2}$}}
\centerline{
  \subfigure[Remove parallel bigons]
  {\label{fig:remove-parallel-bigon}$\mfigb{simplify-20}\to\mfigb{simplify-21}\to\mfigb{simplify-22}\to\mfigb{simplify-23}$}}
\centerline{
  \subfigure[Comb anti-parallel strand]
  {\label{fig:comb-anti-parallel}$\mfigb{simplify-30}\to\mfigb{simplify-31}\to\mfigb{simplify-32}$}}
\centerline{
  \subfigure[Comb parallel strand]
  {\label{fig:comb-parallel}$\mfigb{simplify-40}\to\mfigb{simplify-41}\to\mfigb{simplify-42}$}}
\caption{Steps in the proof of Lemma~\ref{lem:straighten}}
\end{figure}

\begin{lemma}\label{lem:straighten}
  Let $D$ be a triple diagram and let $I$ be a minimal interval of
  the associated pairing.  Then $D$ is related, by a sequence of
  $2\leftrightarrow2$, $1\to0$ and loop-dropping moves, to a diagram $D'$ in which
  there is a strand boundary-parallel along~$I$.
\end{lemma}

\begin{proof}
  We proceed by induction on the number of crossings and number of
  loop components in~$D$, repeatedly applying the lemma to regions
  contained inside of~$D$.  We successively straighten the strand~$S$
  which connects the endpoints of~$I$.
  \begin{description}
  \item[Step 1. Remove self-intersections of~$S$] If $S$ has
    self-intersections, consider an innermost loop~$L$ of~$S$: a loop
    of~$S$ which does not intersect itself, although it may intersect
    other portions of~$S$.  In Figure~\ref{fig:remove-loop}, $L$ is
    shown solid and the rest of $S$ is dashed.  Consider the region
    contained inside $L$, including all triple crossings on the
    boundary of~$L$ except for the self-intersection crossing itself.
    This region has fewer triple crossings than the original
    diagram~$D$ (since we omitted the self-intersection), and the
    short interval between the two points where $L$ intersects the
    boundary of the region is minimal, so by induction we can apply a
    sequence of $2\leftrightarrow2$, $1\to0$ and loop-dropping moves
    to reduce $L$ so that it is
    boundary-parallel inside the region.  A single further $1\to0$ move
    removes the original self-intersection of~$S$.

    Repeat this step until no self-intersections of~$S$ remain.
  \item[Step 2. Remove double intersections with $S$] Consider the
    region~$R$ enclosed by the strand~$S$ and the interval~$I$, not
    including the triple crossings along $S$ itself.  If there is some
    strand that intersects $S$ at least twice, then consider a strand
    segment~$T$ that gives a minimal interval along the boundary
    of~$R$.  Because there are fewer triple crossings inside~$R$ than
    there were in~$D$, we can by induction apply moves until $T$ is
    boundary-parallel to~$R$.  If $T$ and $S$ form an anti-parallel bigon,
    we can apply Lemma~\ref{lem:useful-moves}(\ref{it:move-1}) as in
    Figure~\ref{fig:remove-anti-parallel-bigon} to remove the bigon.  If
    $T$ and $S$ form a parallel bigon, we can similarly apply
    Lemma~\ref{lem:useful-moves}(\ref{it:move-3}) and two $0\to2$ moves
    as in Figure~\ref{fig:remove-parallel-bigon} to remove the bigon.
    
    Repeat this step until there are no strands that intersect~$S$
    twice.  The number of triple intersections along $S$
    strictly decreases at each step, so this terminates.
  \item[Step 3. Comb out triple crossings] At this point, the only
    strands left between $S$ and $I$ come from the boundary within the
    interval~$I$ and cross~$S$.  (Here we use the minimality of~$I$:
    no strands connect~$I$ to itself.)  We wish to make $S$ boundary
    parallel: that is, we want all the strands starting from~$I$ run
    directly to~$S$, without any triple crossings first.  Consider the
    strand~$T$ with the leftmost endpoint along~$I$ which violates
    this condition.  By the hypothesis on $T$ and the minimality
    of~$S$, the interval corresponding to~$T$ is minimal inside the
    region~$R$ between~$S$ and the boundary, omitting the strands to
    the left of~$T$ and the triple crossings along~$S$.  Apply
    induction inside~$R$ to reduce~$T$ until it is boundary parallel.
    If $T$ is oriented anti-parallel to $S$, apply
    Lemma~\ref{lem:useful-moves}(\ref{it:move-2}) to make $T$ run
    straight to $S$, as in Figure~\ref{fig:comb-anti-parallel}; if $T$
    is oriented parallel to $S$, instead apply
    Lemma~\ref{lem:useful-moves}(\ref{it:move-3}), as in
    Figure~\ref{fig:comb-parallel}.

    Repeat this step until every strand runs directly from $I$ to~$S$.
  \item[Step 4. Remove disconnected components] Now $S$ is boundary
    parallel, with the exception of possible disconnected components
    of the diagram between $S$ and the boundary.  For each
    disconnected component, consider a region~$R$ that encloses all
    of the component except for one arc~$T$ on the boundary of the
    component.  The diagram contained inside $R$ has no more triple
    crossings than $D$ and at least one fewer loop, so by induction we
    can reduce it until the unique arc runs boundary parallel on
    the side facing $T$.  This strand forms a simple loop with $T$,
    which we can drop.
    
    Repeat this step until there are no disconnected components
    between $S$ and the boundary, and so $S$ is boundary parallel.
    The total number of loops in the diagram decreases at each step,
    so this terminates.\qedhere
  \end{description}
\end{proof}

\begin{proof}[Proof of Theorem~\ref{thm:domino-flip}]
  Start with two minimal diagram~$D$ and~$D'$.  By repeatedly applying
  Lemma~\ref{lem:straighten} to~$D$, we can obtain a standard
  diagram~$D''$ as in Section~\ref{sec:generating}.
  Since $D$ is minimal, we did not use any $1\to0$ moves in this
  process.  Furthermore, any loop dropping move is necessarily
  preceded by a $1\to0$ move since $D$ is connected and $2\leftrightarrow2$ moves
  preserve connectivity.  Therefore we only used $2\leftrightarrow2$ moves.  In the
  same way~$D'$ can be connected to $D''$ by a sequence of $2\leftrightarrow2$
  moves, and so $D$ and $D'$ can be connected by a sequence of $2\leftrightarrow2$
  moves, as desired.
\end{proof}

\begin{corollary}\label{cor:std-minimal}
  Any standard diagram is minimal.
\end{corollary}

\begin{proof}
  This was proved during the proof of
  Theorem~\ref{thm:domino-flip}.
\end{proof}

\begin{proof}[Proof of Theorem~\ref{thm:reduction}]
  By repeatedly applying Lemma~\ref{lem:straighten}, an arbitrary
  diagram can be reduced to a standard diagram by a sequence of
  $2\leftrightarrow2$, $1\to0$, or loop dropping moves.
\end{proof}

\begin{remark}
  The proof of the connectivity of tilings of a region by
  dominoes~\cite{Thurston90:ConwaysTiling} exploits a \emph{height
    function} defined on the vertices of the tiling which changes by
  $\pm1$ along each edge of the tiling and by $\pm3$ along edges of the
  square grid which are not edges of the tiling.  To show the space of
  domino tilings is connected, transform a diagram to an extremal one
  with minimal height function by repeatedly looking for a vertex
  which is a local maximum of the height function and performing a
  domino flip at that vertex.  The extremal diagrams with respect to
  the height function are also standard diagrams in our terminology.
  Furthermore, general triple diagrams admit a multi-dimensional
  height function on the regions of the diagram which changes by 1 in
  one coordinate when you cross an edge. (See
  \cite[Section~1.2]{HS10:MultidimCube}.)  Standard diagrams are
  extremal with respect to an appropriate projection from the
  multi-dimensional height function to~$\RR$.  Both proofs involve
  connecting an arbitrary tiling to an extremal one.

  There is a corresponding multi-dimensional cube recurrence related
  to rhombus tilings, studied by Henriques and Speyer
  \cite{HS10:MultidimCube}.
\end{remark}

\section{Minimal diagrams}
\label{sec:minimal}

Minimal triple diagrams play a somewhat special role in the
proof above.  In order to get a better understanding of minimal
diagrams, we will prove Theorem~\ref{thm:minimal}

\begin{definition}
  A \emph{badgon} in a general diagram of oriented curves, not
  necessarily a triple diagram, is a monogon, parallel
  bigon, or simple loop.
\end{definition}

Theorem~\ref{thm:minimal} says that a diagram is minimal if and only
if it has no badgons.

\begin{proof}[Proof of Theorem~\ref{thm:minimal}]
  First suppose a diagram~$D$ has a monogon or parallel bigon.  In
  either case, we can straighten the strand as in the proof of
  Lemma~\ref{lem:straighten}, Steps~1 and~2
  (Figures~\ref{fig:remove-loop} and~\ref{fig:remove-parallel-bigon})
  respectively.  At the end of the straightening (and possibly
  earlier), there are one or two $1\to0$ moves, which strictly reduces
  the number of triple crossings; thus~$D$ was not minimal.

  For the other direction, we need a lemma.
  
  \begin{lemma}
    \label{lem:badgon-unchanged}
    Let $D$ and $D'$ be general diagrams, not necessarily connected or
    with only triple crossings, related by one of the following moves:
    \begin{enumerate}\renewcommand\theenumi{\alph{enumi}}
    \item An anti-parallel self-tangency involving distinct strands
      $\mfigb{moves-50}\to\mfigb{moves-51}$;
    \item Perturbing a triple crossing into 3 double crossings
      $\mfigb{moves-60}\to\mfigb{moves-61}$; or
    \item A $2\leftrightarrow2$ move
      $\mfigb{moves-0}\to\mfigb{moves-1}$.
    \end{enumerate}
    Then $D$ and $D'$ either both have no badgons or both have at
    least one badgon.
  \end{lemma}

  \begin{proof} We consider each case in turn.
    \begin{enumerate}\renewcommand\theenumi{\alph{enumi}}
    \item There are two crossings created by the self tangency; call
      them $x$ and $y$.  All other crossings are the same in the two
      diagrams.  Any badgon that does not involve $x$ or $y$ will be
      the same on the two sides of the diagram, so let us suppose
      there is none.  Because the two strands are different, there is
      no monogon created in~$D'$.  Because the self-tangency is
      anti-parallel, we can only have a parallel bigon in~$D'$
      involving both $x$ and $y$ if one of the two strands
      is a loop, which necessarily has badgon in~$D$ as well.  The
      only remaining possibility for a badgon in~$D'$ that is not
      in~$D$ is a parallel bigon involving (wlog) $x$ and another
      crossing elsewhere in the diagram, say $z$, like this:
      \[
      \mfigb{simplify-50}
      \]
      In~$D'$, one strand~$S$ runs from $x$ through $y$ to $z$; the
      other strand~$T$ runs in one direction from $x$ to $z$.  In the
      other direction~$T$ runs from $x$ to $y$.  After the
      intersection with $y$, $T$ is inside the bigon formed by~$S$
      and~$T$ and must exit somewhere.  If it exits by crossing~$T$,
      there is a monogon; if it exits by crossing~$S$, there is
      another parallel bigon.  In either case, this new badgon exists
      in both $D$ and $D'$.
    \item If any two of the strands involved in the triple crossing
      are the same, there is a monogon in both $D$ and $D'$.
      Otherwise, a badgon in~$D'$ will only involve at most one of the three
      crossings and there gives a badgon in~$D$ as well, and vice versa.
    \item If any two of the strands involved are the same,
      there is a badgon on both sides of the diagram, so let us
      suppose all strands are distinct.  The $2\leftrightarrow2$ move can be
      decomposed into a sequence of the moves above: perturb the two
      triple crossings, perform two anti-parallel self-crossings, and
      collect the crossings into two new triple crossings.  None of
      these steps change the presence of badgons.\qedhere
    \end{enumerate}
  \end{proof}
  
  We now finish the proof of Theorem~\ref{thm:minimal}. By
  Theorem~\ref{thm:reduction}, any connected non-minimal
  diagram~$D$ can be turned into a minimal diagram by a sequence of
  moves, necessarily involving a $1\to0$ move.  Let~$D'$ be the
  diagram just before the first $1\to0$ move.  $D'$ has a monogon at the
  site of the $1\to0$ move.  $D$ and~$D'$ are related by a sequence
  of $2\leftrightarrow2$ moves, so by Lemma~\ref{lem:badgon-unchanged}, $D$ has a
  badgon as well.
\end{proof}

\section{Connections and future directions}
\label{sec:connections}

\subsection{Invariants of plane curves and Legendrian knots}
\label{sec:plane-curves}
One natural connection is with Arnol'd's theory of invariants of plane
curves~\cite{Arnold94:TopologicalInvariants}.  Deforming a generic
smooth oriented plane curve in a generic way, you pass through three
types of non-generic points (what Arnol'd calls \emph{perestroikas}):
triple crossing, direct self-tangency, and inverse self-tangency.  If
an invariant of generic plane curves does not change under a triple
crossing, it can be extended naturally to an invariant of curves with
triple points.  If, in addition, it does not change under an inverse
self-tangency, it will not change under the domino flip move.

For instance, Arnol'd's $J_+$ invariant of plane curves, the simplest
invariant unchanged under these two operations, when evaluated on a
curve with only alternating triple points and no double points gives
the number of triple points plus the index of the curve, plus or
minus~1 (depending on the orientation on the boundary of the exterior
region of the plane).

One way to construct invariants of plane curves that do not change
under a triple crossing or an inverse self-tangency is by taking an
invariant of the \emph{Legendrian knot} in $\RR^2 \times S^1$ obtained by
lifting the curve to the unit tangent bundle of the plane.  This
connection with Legendrian knots leads to some natural questions:

\begin{question}
  Which Legendrian knots have planar representations with only
  alternating triple points?
\end{question}

\begin{remark}
  The answer to this question is not ``all'', since the $J_+$
  invariant must be strictly positive for non-trivial knots
  represented by an alternating triple diagram.
\end{remark}

\begin{question}
  Can two different triple diagrams for a single Legendrian knot
  always be related by $2\leftrightarrow2$ moves?
\end{question}

\subsection{Cluster algebras and the octahedron recurrence}
\label{sec:cluster-algebra}

(This section is joint work with Andre Henriques.)

There are two natural cluster algebras as introduced by Fomin and
Zelevinsky~\cite{FZ01:Cluster-Algebras-I} associated to triple
diagrams with a given connectivity.  Loosely speaking, a \emph{cluster
  algebra}~$A$ of rank $n$ is a commutative algebra with unit and no
zero divisors with a distinguished set of \emph{clusters}: sets of
elements $x_1,\ldots,x_n \in A$ so that every element in~$A$ can be written
as a Laurent polynomial in the~$x_i$ in any given cluster.  (A Laurent
polynomial in the $x_i$ is a polynomial in the $x_i$ and $x_i^{-1}$.)
Furthermore, any two clusters can be connected by a sequence of
\emph{exchange relations}: replacing one element~$x_i$ of a cluster by
a conjugate variable~$\tilde x_i$, related by
\[
x_i \tilde x_i = M_1 + M_2
\]
where $M_1$ and $M_2$ are monomials in the~$n-1$ other elements of the
cluster.  (This is a description of cluster algebras, not a
definition.)

Each connectivity class of triple diagrams gives a cluster algebra,
with each triple diagram giving a cluster.  The variables in the
cluster correspond to the white regions in the checkerboard coloring
of the triple crossing diagram.

When you perform a $2\leftrightarrow2$ move, the cluster variables are unchanged if
the region in the center of the move is black.  If the region is
white, the conjugate variable to the variable in the bigon is given by
\begin{equation}\label{eq:cluster-move}
\mfigb{moves-40} \longrightarrow \mfigb{moves-41}
\end{equation}
where
\[
f = \frac{ac + bd}{e}.
\]
One initially surprising fact is the \emph{Laurent phenomenon}: start
from an arbitrary triple diagram and apply a number of these $2\leftrightarrow2$
moves.  Express the variables on each region in terms of the variables
on the original regions using the exchange relation.  Each variable
will be a Laurent polynomial!  The number of terms in these Laurent
polynomials becomes quite large, and you have to divide by earlier
polynomials in the exchange relation, but in each case the division
comes out exactly.  Furthermore, the coefficients of the Laurent
polynomials are positive integers.

This cluster algebra turns out to be a disguised version of the
$n$-dimensional octahedron recurrence~\cite[Section~1.2]{HS10:MultidimCube}.
The coefficients are conjecturally positive, but there is not
currently a combinatorial interpretation in all cases. For the
3-dimensional octahedron recurrence, a combinatorial interpretation in
terms of matchings was given by David
Speyer~\cite{Speyer07:Perfect-Matchings}.

Postnikov~\cite{Postnikov06:TotalPositivity,Scott06:Grassmannians}
has used closely related diagrams to parameterize bases for totally
positive cells in Grassmannians.  Each connectivity corresponds to a
particular totally positive cell; each triple crossing diagram gives a
basis for the Grassmannian.  The bases corresponding to triple
crossing diagrams that are related by a $2\leftrightarrow2$ move are again related
by (\ref{eq:cluster-move}). See
\cite[Section~14]{Postnikov06:TotalPositivity} for a description of
the relation.

\subsection{Relations between relations}
\label{sec:rels-rels}

If Theorem~\ref{thm:generate} guarantees that triple diagrams
generate all possible permutations and Theorem~\ref{thm:reduction}
guarantees that the $2\leftrightarrow2$, $1\to0$, and loop dropping moves give all
relations between triple diagrams, the next natural question
is what the relations between relations are.  For simplicity, we will
restrict ourselves to minimal triple diagrams.

\begin{figure}[htbp]
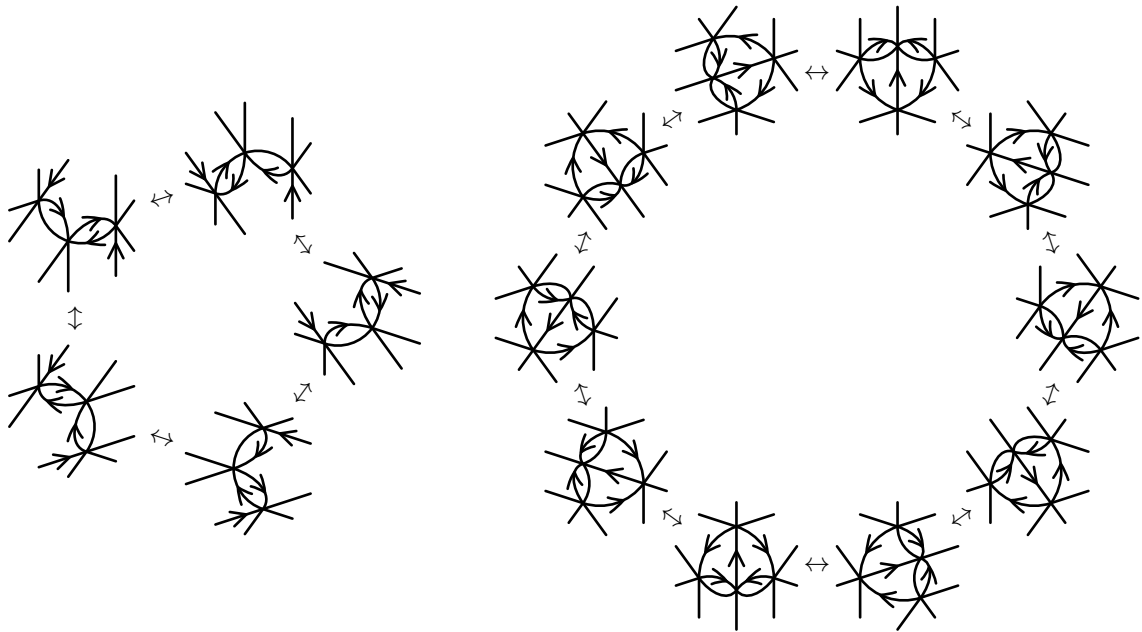

  \centering
  \subfigure{\label{rel-rels-i}$\mathcenter{\includegraphics{crystal-6}}$}
  \qquad
  \subfigure{\label{rel-rels-ii}$\mathcenter{\includegraphics{crystal-16}}$}
  \caption{Conjectured relations between relations}
  \label{fig:rels-rels}
\end{figure}

Before stating a conjecture, we first pick out some particularly
strand connectivities to consider.
\begin{definition}
  Fix an integer $n \ge 1$ and an odd integer $1 \le k \le 2n-1$. The
  \emph{$(n,k)$ strand connectivity} is the connection of strands that
  connects each input to the output position that is $k$ steps
  counterclockwise around the outside of the diagram. (Since $k$ is
  odd, this connects each input strand to an output strand.) Here are
  some examples of connectivities and their minimal diagrams.
  \begin{itemize}
  \item In the connectivities $(n,1)$ and $(n,2n-1)$, each strand is
    connected to a neighbor, and we get a diagram with no triple-crossings.
  \item The connectivity $(3,3)$ gives a single triple-crossing.
  \item In the connectivities $(4,3)$ and $(4,5)$, there are two
    minimal diagrams, each with two triple-crossings. These two
    diagrams are the two sides of the $2 \leftrightarrow 2$ move.
  \item The minimal diagrams for the connectivities $(5,3)$ and
    $(5,5)$ are show in Figure~\ref{fig:rels-rels}.
  \end{itemize}
\end{definition}

\begin{proposition}\label{prop:count-boxes}
  The number of triple-crossings in a minimal diagram of connectivity
  $(n,2\ell+1)$ is $\ell(n-\ell-1)$.
\end{proposition}
\begin{proposition}\label{prop:max-connectivity}
  Among all connectivities for $n$-strand triple-crossing diagrams,
  those with the maximal number of triple-crossings are
  \begin{itemize}
  \item if $n$ is odd: the connectivity $(n,n)$, with $(n-1)^2/4$
    triple-crossings, and
  \item if $n$ is even: the connectivities $(n,n-1)$ and $(n,n+1)$,
    with $n(n-2)/4$ triple-crossings.
  \end{itemize}
\end{proposition}
\begin{proof}[Proof of Propositions~\ref{prop:count-boxes}
  and~\ref{prop:max-connectivity}]
  Straightforward from Theorem~\ref{thm:crossing-count}.
\end{proof}

\begin{conjecture}\label{conj:rels-rels}
  For triple diagrams with a given connectivity, consider the
  2-complex~$C$
  with vertices given by minimal triple diagrams, edges given
  by $2\leftrightarrow2$ moves, and 2-cells of the following types:
  \begin{enumerate}
  \item Squares for two different $2\leftrightarrow2$ moves in non-interfering parts of the
    diagram;
  \item Pentagons as on the left of Figure~\ref{fig:rels-rels} or
    their orientation reverse, from connectivities $(5,3)$ and $(5,7)$; and
  \item Decagons as on the right of Figure~\ref{fig:rels-rels}, from
    connectivity $(5,5)$.
  \end{enumerate}
  Then~$C$ is simply connected.
\end{conjecture}

\begin{figure}
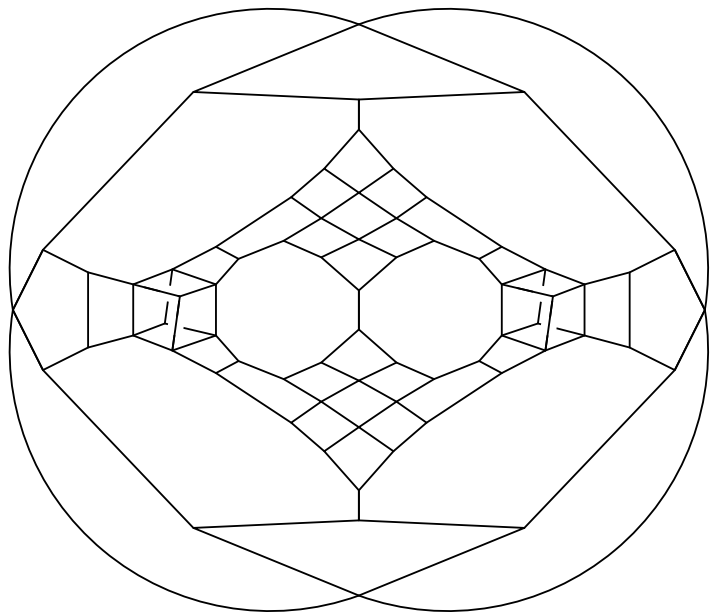

  \[
  \mfig{connectivity-0}
  \]
  \caption{The graph of minimal triple-diagrams for connectivity
    $(6,5)$. Note the presence of two cubes. The graph has 3-fold
    rotational symmetry around these cubes.}
  \label{fig:graph-6-5}
\end{figure}

For instance, the complex $C(n,3)$ for connectivity
$(n,3)$ is the 2-skeleton of the
associahedron $K_n$, which is simply-connected.
The complex $C(6,5)$ is shown in
Figure~\ref{fig:graph-6-5}, and is again simply-connected.

More generally, one could conjecture that there is a contractible
complex with fundamental $n$-cells given by the diagrams of
connectivity $(n+3,k)$ for $k$ an odd number between $3$ and $2n+3$.
As in Conjecture~\ref{conj:rels-rels}, there are also products
of lower-dimensional cells.

\subsection{Planar algebras}
\label{sec:planar-algebras}
The original motivation for considering triple crossings came from
\emph{planar algebras} as introduced by
Jones~\cite{Jones99:Planar-Algebras-I}.  Briefly, planar algebras in
the form we consider have a complex vector space~$B_n$ for each $n$,
the space of \emph{$n$-boxes}.  Each $n$-box conceptually has $2n$
strands attached to the boundary, alternating in and out.  They may be
joined together by oriented \emph{arc diagrams} like this one:
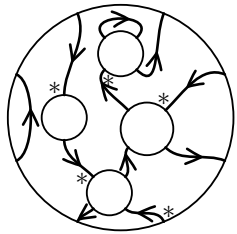
\[
\mfigb{arc-diagram-0} 
\]
This diagram, for instance, has 4 inputs (the little circles), with
4, 2, 4, and 4 strands, respectively; the output (the large circle) has
8 strands.  This diagram therefore gives a multilinear map
\[
B_2 \otimes B_1 \otimes B_2 \otimes B_2 \to B_4.
\]
Each arc diagram gives a similar multilinear map, and if you plug the
output of one multilinear map into the input of another, you get the
multilinear map corresponding to gluing together the two arc diagrams.
Formally, this is the structure of a (typed) operad.  (To be precise,
you must also provide some markings to know the orientations when you
plug in elements.)

Typically one assumes that $B_0$ is 1-dimensional and is identified
with~$\CC$.  In particular, this implies that a simple loop, which can
be thought of as an element of $B_0$, has some definite value $\delta\in\CC$.%
\footnote{One of the subtleties we're ignoring is that a clockwise and
  counterclockwise loop could potentially have different values.}

One natural way to start to classify planar algebras is by their
generating sets.  Every planar algebra contains the
\emph{Temperley-Lieb algebra}, consisting of diagrams with only
non-crossing arcs and no inputs. Non-crossing diagrams with no loops
give us $1, 1, 2, 5, 14, \ldots$ elements inside~$B_n$ for $n=0,1,2,3,4,\ldots$.
Most planar algebras are not spanned by the Temperley-Lieb algebra,
but are generated by the addition of some extra elements.  For
example, Bisch and
Jones~\cite{BJ00:Singly-Generated,BJ03:Singly-Generated-II} have
studied planar algebras that are generated by a single 2-box.
Consider instead planar algebras where the 0-, 1-, and 2-boxes are
spanned by Temperley-Lieb elements and which are generated by a single
3-box.  If we draw the 3-box as a triple crossing,%
\footnote{Another subtlety we're ignoring here is that the 3-box need
  not be symmetric under rotation by 120 degrees, as the triple
  crossing is.  This turns out not to affect the subsequent
  discussion.}  then every triple diagram gives an $n$-box.  If we
assume further that the dimension of~$B_4$ is less than or equal to 25
and that the relations are ``generic'' and ``symmetric'' in senses to
be made precise below, the results in this paper imply that the
dimension of~$B_n$ is less than or equal to $n!$ for all~$n$.

To prove this, consider diagrams modulo ``simpler'' diagrams that
either have fewer loops or fewer triple crossings, and look at the
three moves used in Theorem~\ref{thm:reduction}.
\begin{itemize}
\item A loop is equal to a constant~$\delta$.
\item Because $B_2$ is spanned by Temperley-Lieb elements, the
  diagram~$\mfigb{moves-11}$ is equal to a linear
  combination of diagrams with no crossings, giving the $1\to0$ move
  modulo lower order terms.
\item There are 26 elements of $B_4$ obtained from the 14 elements of
  the Temperley-Lieb algebra, the 8 ways of adding a single arc to a
  single triple crossing, and the 4 ways of connecting two triple
  crossings along two adjacent arcs, so by assumption there is at
  least a 2-dimensional space of relations between these 26 elements.
  The highest terms are of the last kind.  Assume that the relations
  are generic, in the sense that when we take the quotient by the
  space spanned by the
  first 22 diagrams, the space of relations is a nontrivial
  subspace~$W$ of the 4-dimensional space~$V$ with basis given by the
  remaining 4 diagrams.

  Rotating by 90 degrees gives a symmetry of
  order two of~$V$ that must preserve~$W$.  $V$ decomposes into two
  subspaces~$V = V^+ \oplus V^-$, both of dimension two.  If we furthermore
  assume that~$W$ is symmetric when we reverse all the arrows, $W$
  must contain either~$V^+$ or~$V^-$.  Then we have either
  \begin{equation}\label{eq:rel-plus}
  \mfigb{moves-0} = \mfigb{moves-1} + (\text{lower order terms})
  \end{equation}
  or
  \begin{equation}\label{eq:rel-minus}
  \mfigb{moves-0} = -\mfigb{moves-1} + (\text{lower order terms}),
  \end{equation}
  respectively.
\end{itemize}
The loop in Figure~\ref{rel-rels-i} shows that Equation~\eqref{eq:rel-minus}
implies that the entire space of diagrams of two triple crossings
vanishes modulo lower order terms. Thus the three moves involved in
Theorem~\ref{thm:reduction} are all true modulo lower
order terms.  Therefore, modulo lower order terms,
\begin{itemize}
\item Non-minimal diagrams are equal to 0,
  and
\item Minimal diagrams with the same connectivity are equal to each
  other.
\end{itemize}
(We need the fact that you don't need to increase the size of the
diagram in Theorem~\ref{thm:reduction}.)
Therefore $B_n$ is spanned by $n!$ minimal diagrams, one for each
connectivity.

One source of planar algebras satisfying these hypotheses is the skein
module of the HOMFLYPT polynomial, a well-known invariant of oriented
links.  You might expect a knot invariant to give a planar algebra
based on a simple crossing, which would be a 2-box; however, the
orientations in the HOMFLYPT polynomial and the requirement that the
strands alternate in-out force you to consider triple crossings.

Since the first version of this paper was published, this approach has
been worked out by Jones, Liu, and Ren \cite{JLR16:ThurstonRel}.

\bibliographystyle{hamsalpha}
%\bibliography{topo,cluster,geom,planaralg,curves,graphs}
\bibliography{Dominoes}
\end{document}